# Exact arithmetic as a tool for convergence assessment of the IRM-CG method

Josip Dvornik[1]    Antonia Jaguljnjak Lazarevic[2]    Damir Lazarevic[1,3]    Mario Uroš[1]

[1]University of Zagreb, Faculty of Civil Engineering, Kaciceva 26, Zagreb 10 000, Croatia

[2]University of Zagreb, Faculty of Mining, Geology and Petroleum Engineering, Pierottijeva 6, Zagreb 10 000, Croatia

**Abstract**

Using exact computer arithmetic, it is possible to determine the (exact) solution of a numerical model without rounding error. For such purposes, a corresponding system of equations should be exactly defined, either directly or by rationalisation of numerically given input data. In the latter case there is an initial roundoff error, but this does not propagate during the solution process. If this system is first exactly solved, then by the floating-point arithmetic, convergence of the numerical method is easily followed. As one example, IRM–CG, a special case of the more general Iterated Ritz method and interesting replacement for a standard or preconditioned CG, is verified. Further, because the computer demands and execution time grow enourmously with the number of unknowns using this strategy, the possibilities for larger systems are also provided.

*Keywords*: exact arithmetic, benchmark, rounding error, Iterated Ritz method, Conjugate gradient method

## 1 Introduction: The Iterated Ritz Method

The Iterated Ritz Method (IRM) is an iterative approach to the solve of the symmetric positive definite (SPD) system $\mathbf{A}\mathbf{x} = \mathbf{b}$ based on successive minimisation of the corresponding energy (the quadratic function)

$$f(\mathbf{x}) = \frac{1}{2}\mathbf{x}^{\mathrm{T}}\mathbf{A}\mathbf{x} - \mathbf{x}^{\mathrm{T}}\mathbf{b} \tag{1}$$

inside a small subspace formed at each step [1]. The main strategy is to present solution increments by the Ritz idea:

$$\mathbf{p}_{(i)} = \mathbf{\Phi}_{(i)}\mathbf{a}_{(i)} \tag{2}$$

where $\mathbf{\Phi}_{(i)} = [\,\boldsymbol{\phi}_{1,(i)}\ \boldsymbol{\phi}_{2,(i)}\ \cdots\ \boldsymbol{\phi}_{m,(i)}\,]$ is a matrix of linearly independent coordinate vectors (that form subspace) and $\mathbf{a}_{(i)}$ is the vector of corresponding coefficients. Using this approach, the energy decrement achieved by (2), is also presented in the quadratic form

$$\Delta f(\mathbf{a}_{(i)}) = \frac{1}{2}\mathbf{a}_{(i)}^{\mathrm{T}}\overline{\mathbf{A}}_{(i)}\mathbf{a}_{(i)} - \mathbf{a}_{(i)}^{\mathrm{T}}\bar{\mathbf{r}}_{(i)} \tag{3}$$

where $\overline{\mathbf{A}}_{(i)} = \mathbf{\Phi}_{(i)}^{\mathrm{T}}\mathbf{A}\mathbf{\Phi}_{(i)}$ and $\bar{\mathbf{r}}_{(i)} = \mathbf{\Phi}_{(i)}^{\mathrm{T}}\mathbf{r}_{(i)}$ are the SPD generalised (Ritz) matrix and residual vector, both of order $m$. Minimisation of (3) leads to a system of equations that should be solved at each step:

$$\overline{\mathbf{A}}_{(i)}\mathbf{a}_{(i)} = \bar{\mathbf{r}}_{(i)} \tag{4}$$

This is a very small system, because only several coordinate vectors are applied ($m \ll n$). The solution is used to find the increment in (2), and $\mathbf{x}_{(i+1)} = \mathbf{x}_{(i)} + \omega\mathbf{p}_{(i)}$ is subsequently updated.

Coefficent $\omega \in (0, 2)$ is the relaxation factor, known from the Successive Overrelaxation method, and might improve the convergence. Some optimal $\omega$ exists, but they vary for every unknown and solution step. Moreover, the determination is usually more 'expensive' than the benefit of possible improvement. Generally, $\omega$ may be guessed by intuition or experience, and kept constant during the solution process.

---

[3]   Corresponding author: e-mail: `damir@grad.hr`, 0000-0002-7439-719X, tel.: +385 1 4639 640

The residual is recursively defined as $\mathbf{r}_{(i+1)} = \mathbf{r}_{(i)} - \omega \mathbf{A} \mathbf{p}_{(i)}$ and should always be corrected after some (say $k$) number of steps using the equilibrium relation $\mathbf{r}_{(i+1)} = \mathbf{b} - \mathbf{A}\mathbf{x}_{(i+1)}$, because of accumulated round-off errors. The process is terminated after the convergence criterion is reached, i.e. $\|\mathbf{r}_{(i)}\|_2 \leq \varepsilon \|\mathbf{r}_{(0)}\|_2$, where $\varepsilon$ is a very small positive number. To avoid useless calculations (if the algorithm has trouble converging) the maximum number of steps ($n_{\max}$) should also be defined. Simple pseudocode, with sequence of instructions common for iterative solution methods, is given by the Algorithm 1.

Therefore, at each step, coordinate vectors spanning the subspace are created, within which the energy of the system is reduced. This is why the small system (4) needs to be solved (most often by some direct solver). If $\omega = 1$ it is the largest reduction (local energy minimum), which is not necessarily optimal for global convergence. The procedure leads to a gradient class solution method that combines iterative and direct solution strategies. If the iterative process is convergent, the sum of small-system solutions approaches the large (original) system solution, and the sum of small-system energies monotonically decreases and approaches the minimum of the large system [2].

For the convergence of IRM one coordinate vector not orthogonal to the current residual is sufficient. It can be $\mathbf{r}_{(i+1)}$ itself, or multiplied by some SPD matrix. A previous solution increment $\mathbf{p}_{(i)}$ contributes to faster convergence, and it is also frequently used. This vector is known from the previous step, therefore it 'costs' nothing, contrary to other vectors that should be generated somehow. Such vectors are of different efficiency, but are very freely selected – they should only be linearly independent.

IRM can also be considered as a generalisation of some iterative methods [3]. Depending on the choice of coordinate vectors, solvers can be represented or interpreted as special cases of this approach. Furthermore, it is possible to combine good properties of several methods simultaneously. If appropriate vectors are selected, convergence should proceed faster than using any single method considered. Here, an improved CG algorithm (IRM-CG) is briefly presented [4], which is applicable to sparse and large SPD systems arising in many physics applications.

## 2　Non-recursive CG-like algorithm without the need to restart

IRM-CG also starts with the Steepest Descent step. Other steps are executed using a simple alternative to CG simulated by IRM with two coordinate vectors: $\mathbf{r}_{(i+1)}$ and $\mathbf{p}_{(i)}$. Vectors span a two-dimensional subspace. At each step, a system of two equations is solved and a local energy minimum within that plane is found (Algorithm 2). As in CG, only one matrix-vector multiplication is required per step (line 13, Algorithm 2), with appropriate transformations.

Local energy minima are numerically 'exact', contrary to standard CG where the solution of two equations is sought by the equivalent recursive $\mathbf{A}$-orthogonalisation. Because of the accumulation of round-off errors, orthogonality only really exists for a few adjacent vectors and convergence difficulties are present for large ill-conditioned systems. Many restart [5] and preconditioning techniques [6–10] improve convergence.

In IRM-CG, error in $\mathbf{A}$-orthogonality also exists, but orthogonality is not used nor accumulated during an iterative process. Therefore, inherited errors decrease, though non-exact arithmetic (as

**Algorithm 1** Basic IRM algorithm

**Require:** $\mathbf{A}$, $\mathbf{b}$, $\mathbf{x}_{(0)}$, $\omega$, $\varepsilon$, $k$, $n_{max}$   {usually $\mathbf{x}_{(0)} \leftarrow \mathbf{0}$}

**Ensure:** $\mathbf{x}_{(i+1)}$   {close to $\mathbf{x}$}

1:  $i \leftarrow 0$   {initialisation: steepest descent}

2:  $\mathbf{r}_{(0)} \leftarrow \mathbf{b} - \mathbf{A}\mathbf{x}_{(0)}$

3:  $q \leftarrow \mathbf{r}_{(0)}^{\mathrm{T}}\mathbf{r}_{(0)} / (\mathbf{r}_{(0)}^{\mathrm{T}}\mathbf{A}\mathbf{r}_{(0)})$

4:  $\mathbf{p}_{(0)} \leftarrow q\mathbf{r}_{(0)}$

5:  **while** $(\|\mathbf{r}_{(i)}\|_2 > \varepsilon\|\mathbf{r}_{(0)}\|_2) \wedge (i \leq n_{\max})$ **do** {iterated Ritz method}

6:    $\mathbf{x}_{(i+1)} \leftarrow \mathbf{x}_{(i)} + \omega\mathbf{p}_{(i)}$

7:    **if** $i \bmod k \neq 0$ **then**

8:      $\mathbf{r}_{(i+1)} \leftarrow \mathbf{r}_{(i)} - \omega\mathbf{A}\mathbf{p}_{(i)}$

9:    **else**

10:     $\mathbf{r}_{(i+1)} \leftarrow \mathbf{b} - \mathbf{A}\mathbf{x}_{(i+1)}$

11:   **end if**

12:   generate $[\boldsymbol{\phi}_{1,(i)}\ \boldsymbol{\phi}_{2,(i)}\ \cdots\ \boldsymbol{\phi}_{m,(i)}]$

13:   $\overline{\mathbf{A}}_{(i)} \leftarrow [\boldsymbol{\phi}_{1,(i)}\ \boldsymbol{\phi}_{2,(i)}\ \cdots\ \boldsymbol{\phi}_{m,(i)}]^{\mathrm{T}}\mathbf{A}[\boldsymbol{\phi}_{1,(i)}\ \boldsymbol{\phi}_{2,(i)}\ \cdots\ \boldsymbol{\phi}_{m,(i)}]$

14:   $\overline{\mathbf{r}}_{(i)} \leftarrow [\boldsymbol{\phi}_{1,(i)}\ \boldsymbol{\phi}_{2,(i)}\ \cdots\ \boldsymbol{\phi}_{m,(i)}]^{\mathrm{T}}\mathbf{r}_{(i+1)}$

15:   $\mathbf{a}_{(i)} \leftarrow \overline{\mathbf{A}}_{(i)}^{-1}\overline{\mathbf{r}}_{(i)}$

16:   $\mathbf{p}_{(i+1)} \leftarrow [\boldsymbol{\phi}_{1,(i)}\ \boldsymbol{\phi}_{2,(i)}\ \cdots\ \boldsymbol{\phi}_{m,(i)}]\mathbf{a}_{(i)}$

17:   $i \leftarrow i + 1$

18: **end while**   {end iterated Ritz method}

in every numerical process) affects (but does not threaten) the convergence. As a consequence, restarting of IRM-CG is not needed. Further, preconditioning-like techniques [11] can be adopted easily [4].

One more advantage of this formulation is natural adoption of the relaxation factor. Using $\omega \neq 1$, $\mathbf{A}$-orthogonality is lost, which contradicts the standard CG algorithm, which is just based on the $\mathbf{A}$-orthogonality. Therefore, only preconditioners of classical CG could be modified by $\omega$.

## 3   IRM-CG is equivalent to CG

These methods are equivalent, and it is possible for them to be interchanged. Each step may be performed by CG or IRM-CG, regardless of how the earlier steps were realised. If CG is preferred, we suggest that a single IRM-CG step is ocassionally executed, before the orthogonality error becomes too large. This could be termed 'refresh' rather than 'restart'.

If exact arithmetic is considered, IRM-CG and CG have an identical sequence of step results. The exact solution is obtained after the total number of $m$ steps, where $m$ is the number of different 'active' eigenvalues [12]. If $\mathbf{b}$ is represented as a sum of eigenvectors $\mathbf{v}_j$, i.e. $\mathbf{b} = \sum a_j \mathbf{v}_j$, eigenvectors (and corresponding eigenvalues) with $a_j \neq 0$ are 'active' ('inactive' otherwise). Of course, $m$ can be found only if all $n$ eigenpairs are detected. Multiple eigenvalues should be counted as one, and 'inactive' eigenvalues are not counted at all. This comment is of less practical significance, because IRM-CG is interesting as an iterative, not as a direct solution method [13].

## 4   Simple illustrative example

The above considerations are easily proved by the exact arithmetic. Consider a simple example defined exactly by the direct stiffness method (Figure 1). Elements of $\mathbf{A}$ and $\mathbf{b}$ are rational numbers and integers. A system of equations is solved by CG and IRM-CG with the exact arithmetic (subscript E). Thus, rounding error is completely avoided. Both algorithms are realised with the Wolfram Language [14] used in the Wolfram Mathematica, version 11.3 [15]. Results at every step, (such as residuals, displacements, reactions or internal forces) are also exact and are the same by both methods. The system has 8 DoF and 6 active eigenvectors, because $a_3$ and $a_7$ are zero (as an integer). Therefore, six steps are needed to obtain the solution. After initialisation and $\|\mathbf{r}_{(0)}\|_2 = 1$, remaining relative residual norms are:

$$\frac{1}{2}, \quad 3\,\frac{\sqrt{15\,881}}{814}, \quad 48\,\frac{\sqrt{1\,134\,618\,045}}{20\,474\,189}, \quad 24\,\frac{\sqrt{\cdots}}{\cdots}, \quad 2\,306\,892\,210\,264\,\frac{\sqrt{\cdots}}{\cdots}, \quad 0 \qquad (5)$$

The fourth and the fifth norm contain very large integers (marked with dots) and are not included here. The sixth norm is exact zero. Final displacement components are

$$\mathbf{x}_{(6)} = \frac{1}{a}\begin{bmatrix} 1\,440 & 11\,440\,004\,167 & 3\,840 & 0 & 22\,880\,008\,334 & -1\,440 & \frac{\cdots}{64} & \cdots \end{bmatrix}^{\mathrm{T}} \qquad (6)$$

where $a = 39\,440\,013\,077$. The last two values are also too large to show them here. It is worth noting that, because of the beam symmetry condition, the fourth DoF (rotation of the second joint) is also the integer zero. Even for a such small system precision deteriorates, if a

---

**Algorithm 2** Basic IRM-CG algorithm

---

**Require:** $\mathbf{A}$, $\mathbf{b}$, $\mathbf{x}_{(0)}$, $\omega$, $\varepsilon$, $k$, $n_{max}$    {usually $\mathbf{x}_{(0)} \leftarrow \mathbf{0}$}

**Ensure:** $\mathbf{x}_{(i+1)}$    {close to $\mathbf{x}$}

1: $i \leftarrow 0$    {initialisation: steepest descent}

2: $\mathbf{r}_{(0)} \leftarrow \mathbf{b} - \mathbf{A}\mathbf{x}_{(0)}$

3: $q \leftarrow \mathbf{r}_{(0)}^{\mathrm{T}}\mathbf{r}_{(0)} / (\mathbf{r}_{(0)}^{\mathrm{T}}\mathbf{A}\mathbf{r}_{(0)})$

4: $\mathbf{p}_{(0)} \leftarrow q\mathbf{r}_{(0)}$

5: $\boldsymbol{\beta}_{(0)} \leftarrow \mathbf{A}\mathbf{p}_{(0)}$

6: **while** $(\|\mathbf{r}_{(i)}\|_2 > \varepsilon\|\mathbf{r}_{(0)}\|_2) \wedge (i \leq n_{\max})$ **do** {IRM-CG method}

7:    $\mathbf{x}_{(i+1)} \leftarrow \mathbf{x}_{(i)} + \omega\mathbf{p}_{(i)}$

8:    **if** $i \bmod k \neq 0$ **then**

9:       $\mathbf{r}_{(i+1)} \leftarrow \mathbf{r}_{(i)} - \omega\boldsymbol{\beta}_{(i)}$

10:    **else**

11:       $\mathbf{r}_{(i+1)} \leftarrow \mathbf{b} - \mathbf{A}\mathbf{x}_{(i+1)}$

12:    **end if**

13:    $\boldsymbol{\alpha}_{(i)} \leftarrow \mathbf{A}\mathbf{r}_{(i+1)}$    {sole matrix-vector multiplication}

14:    $\overline{\mathbf{A}}_{(i)} \leftarrow [\,\mathbf{r}_{(i+1)}\ \mathbf{p}_{(i)}\,]^{\mathrm{T}}[\,\boldsymbol{\alpha}_{(i)}\ \boldsymbol{\beta}_{(i)}\,]$    {$\overline{\mathbf{A}}_{(i)}$ is symmetric: $\mathbf{r}_{(i+1)}^{\mathrm{T}}\boldsymbol{\beta}_{(i)} = \mathbf{p}_{(i)}^{\mathrm{T}}\boldsymbol{\alpha}_{(i)}$}

15:    $\overline{\mathbf{r}}_{(i)} \leftarrow [\,\mathbf{r}_{(i+1)}^{\mathrm{T}}\mathbf{r}_{(i+1)}\ \ \omega\mathbf{r}_{(i+1)}^{\mathrm{T}}\mathbf{p}_{(i)}\,]^{\mathrm{T}}$    {if $\omega = 1$ second term is zero}

16:    $\mathbf{a}_{(i)} \leftarrow \overline{\mathbf{A}}_{(i)}^{-1}\overline{\mathbf{r}}_{(i)}$

17:    $\mathbf{p}_{(i+1)} \leftarrow [\,\mathbf{r}_{(i+1)}\ \mathbf{p}_{(i)}\,]\mathbf{a}_{(i)}$

18:    $\boldsymbol{\beta}_{(i+1)} \leftarrow [\,\boldsymbol{\alpha}_{(i)}\ \boldsymbol{\beta}_{(i)}\,]\mathbf{a}_{(i)}$

19:    $i \leftarrow i + 1$

20: **end while**    {end IRM-CG method}

---

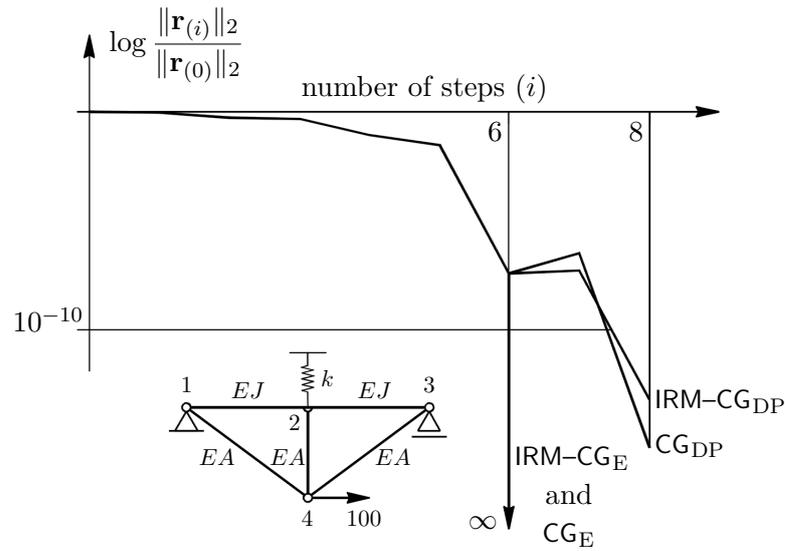

Figure 1: Relative residual norms of the exact (E) and double-precision (DP) implementation of IRM–CG and CG applied to simple structural system

double-precision arithmetic is used (subscript DP). Therefore, the solution is found after two more (eight) steps, providing $\varepsilon \leq 10^{-10}$ is satisified (see Figure 1).

## 5 Check of the algorithm stability

For such a small example, an exact check of the algorithm stability is possible. Because the number of operations is small and the rounding error is negligible, both methods are intentionally perturbed by $\delta = 1$ [4], added (for example) to the seventh component of $\mathbf{p}_{(1)}$ (at the end of the first step):

$$p_{7,(1)} \leftarrow p_{7,(1)} + \delta \tag{7}$$

This effect is similar to the loss of orthogonality, which is common for the iterative methods, if applied to large and ill-conditioned systems. Then, exact arithmetic is again used to obtain the solution. Even with the perturbation, IRM–CG gives an exact result, because $\delta$ is in the loading (residual) direction and therefore lies in the plane spanned by the coordinate vectors. In this particular case, the behaviour of IRM-CG is as if $\delta = 0$ (Figure 2).

Roughly, if $\delta$ is split into two components at each step, the one inside, and the other orthogonal to the plane, the first component is exactly resolved and does not produce inherited error. Using CG, both components causes propagation of error and two more steps are needed to obtain the solution. Similar behaviour is noticed in the double-precision environment (added to the Figure 2). Such perturbations may be 'induced' by the program code of any solution method. Convergence is then verified by comparing the results obtained with exact and floating point implementation.

## 6 More general examples

Consider a larger model – minimally supported (externally statically determinate) cube, loaded with the unit force at the top (Figure 3a). The cube is discretised by a single Lagrangian $\mathcal{C}^0$ finite element with 192 DoF. Inner nodes are not statically condensed and the stiffness matrix is resolved exactly, using rational numbers [16]. Corner stiffnesses are defined similarly and are used to control the condition number $\kappa(\mathbf{A})$, calculated as the ratio of the extreme eigenvalues.

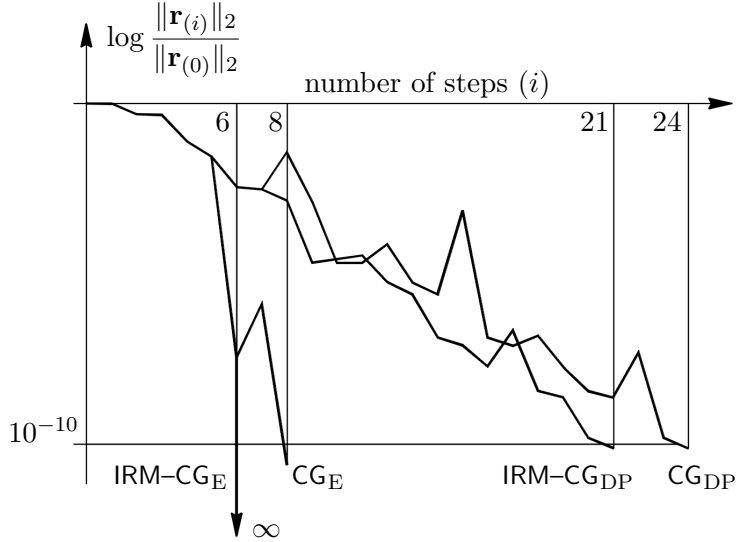

Figure 2: Previous example solved similarly, but with the perturbation of $p_{7,(1)}$

Regardless of $\kappa(\mathbf{A})$, exact relative residual norms are the same for both methods, and the processes stop after $m = 188$ steps. This is because four inactive eigenvalues (eigenvectors orthogonal to the loading) are detected. The last residual norm is integer zero, $\|\mathbf{r}_{(188)}\|_2 = 0$, which means that exact solution of a numerical model is obtained, $\mathbf{x}_{(188)} = \mathbf{x}$. Just for curiosity's sake, vertical deflection of the loaded point (last, 192th component of vector $\mathbf{x}$) is:

$$x_{192,(188)} = -\frac{64464038682410325066453295792043745775932541 3}{288744303487566466737409740217179460197860000}$$

In the numerical environment, methods respond similarly to a well-conditioned model (Figure 3a), but for the ill-conditioned case IRM–CG is more stable, especially if higher accuracy ($\varepsilon = 10^{-10}$) is needed (Figure 3b). Of course, approximate results (for example $\mathbf{x}_{(405)}$ and $\mathbf{x}_{(420)}$, or $\mathbf{x}_{(659)}$ and $\mathbf{x}_{(847)}$) are mutually close and match the exact solution $\mathbf{x}_{(188)}$ reasonably well.

This strategy can also be applied to more complex models (not only) from structural engineering practice, in a combination of various finite elements. If a stiffness matrix is not exactly defined, it is always possible to be rationalized. Elements that are very close to zero could be replaced by the exact integer zero. Using this strategy, the initial roundoff error remains, but it is not accumulated during the solution process (providing exact arithmetic is used). The result is very close to the exact solution of a numerical model $\mathbf{x}_{(m)} \approx \mathbf{x}$. Again, detailed algorithm performance (various steps and final results) can be compared with that obtained by the floating point arithmetic, executed with various numbers of significant digits.

For example, the model from the Figure 4 comprises beam and thin shell elements with displacements and rotations as unknowns. The system has pinned supports and is loaded by two unit moments applied at the centres of the plates. Here, the condition number is controlled by changing the stiffnesses of two corner columns denoted by $\mathsf{A}$. The system has 183 unknowns and all eigenvectors are active (all $a_j \neq 0$, thus $m = n$). Therefore, using exact arithmetic, the solution is found after 183 steps (Figures 4a and 4b). As in previous example, IRM–CG is better for the non-well-posed problem, if higher accuracy has to be achieved.

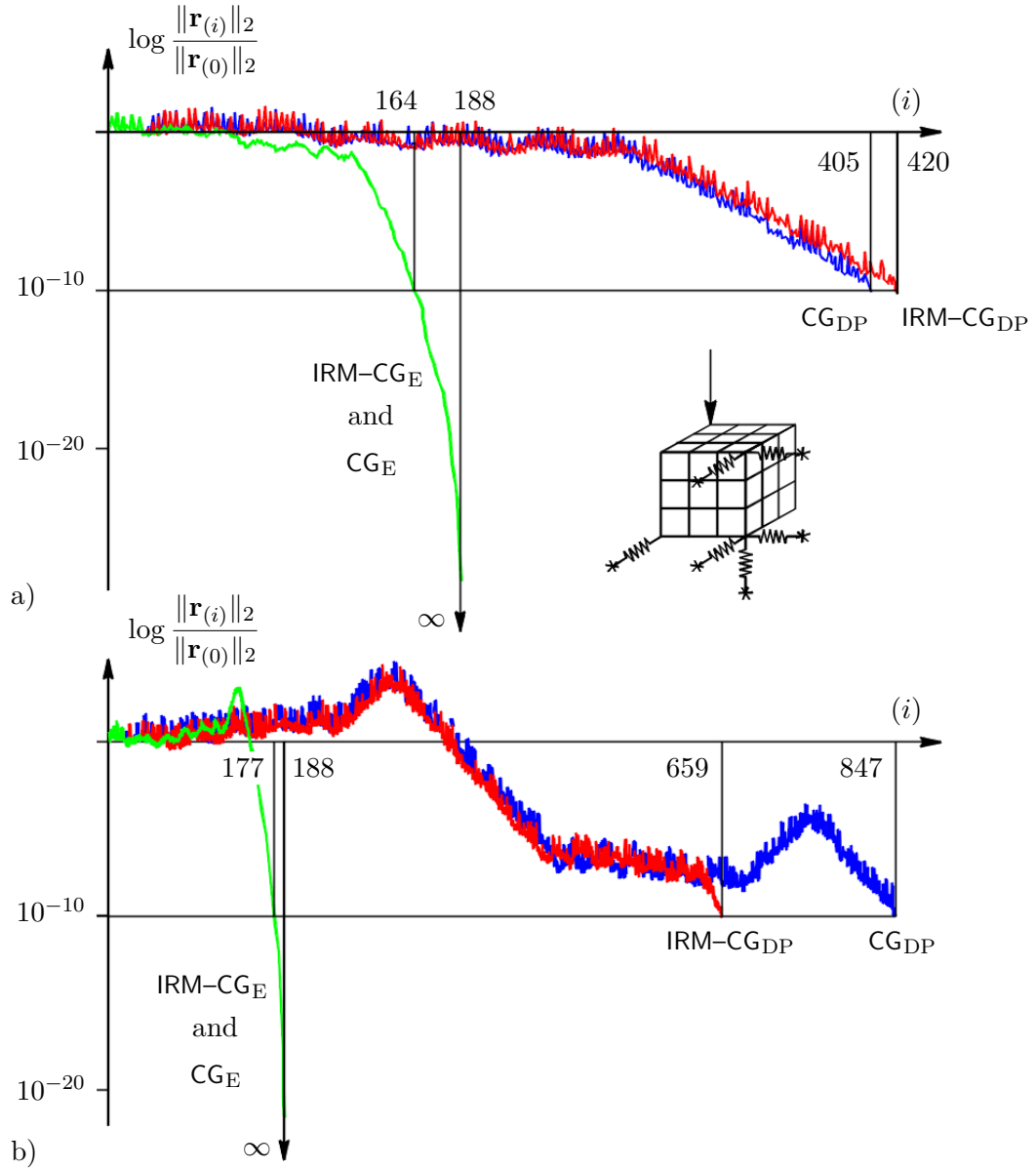

Figure 3: Relative residual norms of the exact and double-precision implementation of CG and IRM–CG applied to FEA of cube model: a) $\kappa(\mathbf{A}) = 7,7 \cdot 10^5$, b) $\kappa(\mathbf{A}) = 6,4 \cdot 10^{12}$

Interestingly, there is a large difference in the number of steps between the exact and double-precision approach, even for a well-conditioned systems. Curves are mutually close only at the early stage of calculation (at the very beginning they practically collide), because rounding error is not accumulated enough.

Obviously, an increase of the accuracy to more than 16 decimal digits of mantissa may be justified. With more significant figures, residual curves that correspond to the numerical implementation of methods are closer to each other, and converge to the curve of the exact

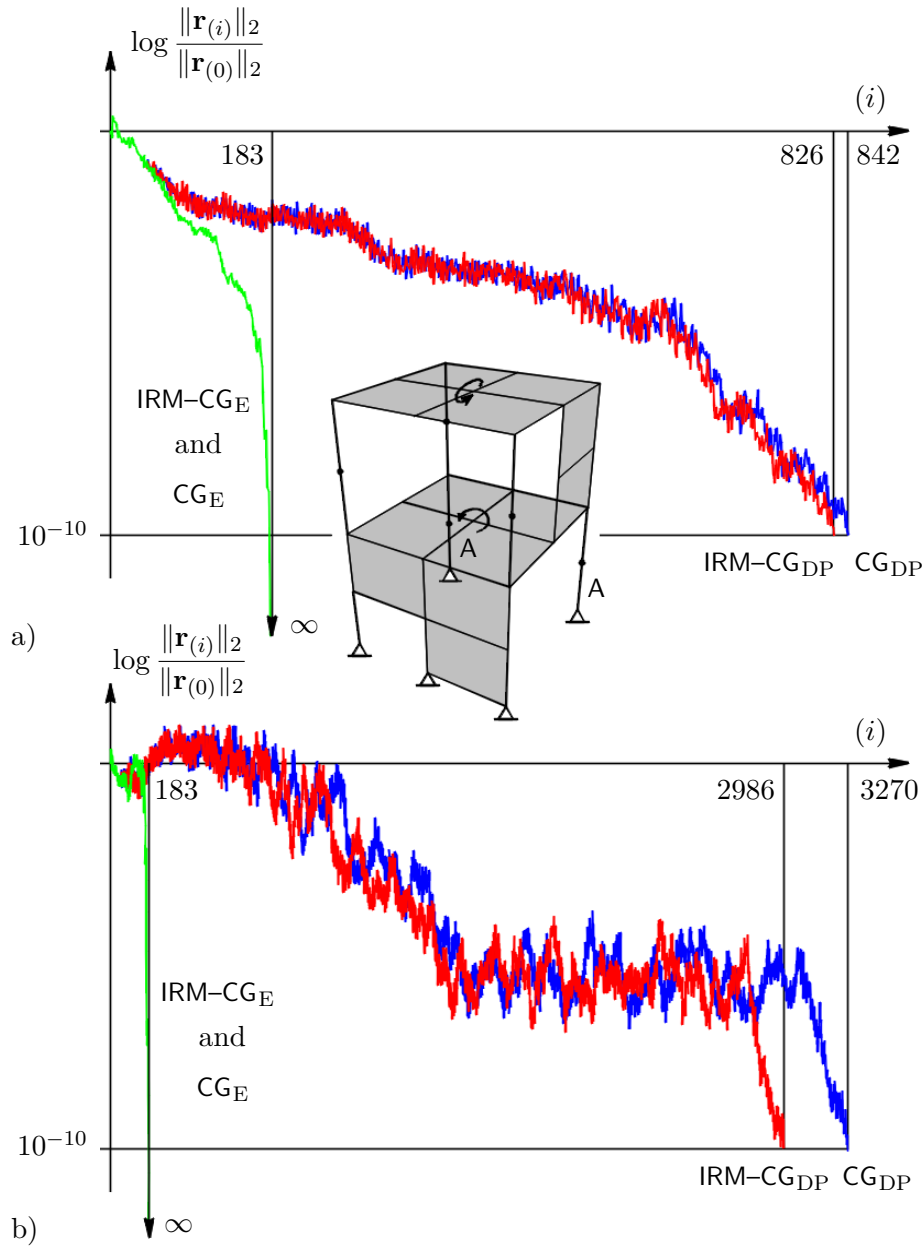

Figure 4: Relative residual norms of the exact and double-precision implementation of CG and IRM–CG applied to FEA of two-story structure: a) $\kappa(\mathbf{A}) = 1,1 \cdot 10^5$, b) $\kappa(\mathbf{A}) = 7,9 \cdot 10^8$

approach. In the limiting case (theoretically, for an infinite number of digits), all three curves must collide. For reasons of clarity, residual functions obtained by the higher precision arithmetic are not added to the figures.

It should be emphasized that, albeit more realistic, a small systems are analysed herein. If the results from the Figures 3 or 4 were to be valid for a large system, the methods would be inefficient. In practice, it is just the opposite – the number of steps is much smaller than the number of unknowns. An explanation is given in Figure 5, where a typical decrease of the residual norm for a large number of DoFs (active eigenvalues $m$) is sketched.

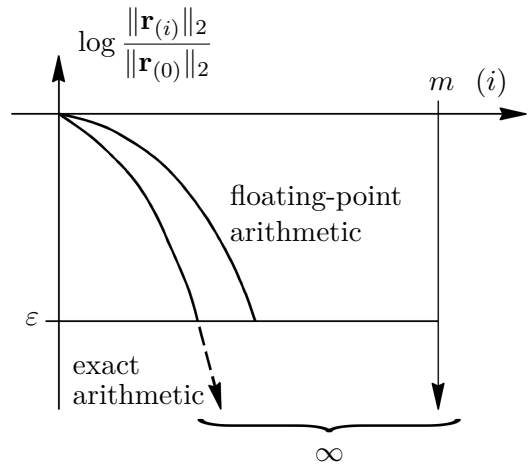

Figure 5: Typical decrease of the relative residual norms for large system

There is a possibility that the curve corresponding to the floating-point arithmetic lies below the curve of the exact approach. In such a scenario, round-off errors are not highly accumulated; they rather very frequently cancel each other during the solution process. However, this is more often the exception rather than the rule.

# 7   Possibilities for large systems

If a large system of equations is considered, more steps, computer memory, and time are required to obtain (not only) an exact solution. Symbolic expressions become very complicated, while whole numbers and fractions grows enormously. Such problems can be overcome by three strategies.

First, instead of solving exactly, arbitrary-precision arithmetic (with more significant figures than fixed-precision arithmetic) may be applied. In that case, not exact, but very accurate solutions may be found.

Second, an inverse method can be used to establish a benchmark. Simply put, a system matrix is multiplied by some solution to find the corresponding right-hand side vector. The input data and matrix-vector multiplication need to be exact. Such an example is then solved by some numerical (here iterative) method, and the results should be easily compared.

Third, systems with a diagonal matrix may provide an interesting approach. This strategy is based on the consequences of the spectral theorem. Briefly, every SPD matrix can be diagonalised as $\mathbf{V}^\mathrm{T}\mathbf{A}\mathbf{V}$, with diagonal elements as eigenvalues and $\mathbf{V} = [\,\mathbf{v}_1 \quad \mathbf{v}_2 \ \ldots \ \mathbf{v}_n\,]$. The spectra of the original and diagonal matrix are the same. The corresponding right-hand side vector is $\mathbf{V}^\mathrm{T}\mathbf{b}$. This transformation is actually a rotation such that the eigenvectors become parallel to the space coordinate axes. If the original system is solved using IRM–CG or CG and steps are then transformed (rotated to that specific position), the results are equal as if the methods had been directly applied to a diagonal system. Therefore, work with a diagonal or original matrix is equivalent (except for the rounding errors if floating-point arithmetic is used).

It should be mentioned that this transformation is based on the eigensolution, which is more 'expensive' than the solution of the corresponding system. Therefore, such a strategy (mostly) makes no sense for real-life engineering problems; it could only be used for clever benchmarks and for tests used to analyse various numerical methods.

Furthermore, for such experiments (with either exact or floating-point arithmetic) the above mentioned transformations of the original matrix are not needed. Put simply, an appropriate diagonal matrix is formed directly. Here 'appropriate' means that the matrix spectrum, condition number, or eigenvalue distribution is easily controlled and that their effects on solver performance could be clearly recognised. Moreover, even for large systems, memory and time requirements are not very demanding and change of the residual spectrum [2, 3] over steps is easy to follow. Of course, purpose of this approach is not to efficiently solve a system, but to get insight into the behaviour of iterative methods. Solution to a diagonal system is easily obtained directly, because equations are mutually independent.

## 8  Conclusion

With the rapid development of computer algebra systems, exact arithmetic has become a very useful tool for convergence assessment of numerical solution methods. Herein, a simple alternative to CG named IRM–CG was tested. The method is not based on conjugacy, therefore it has several advantages over classical CG. It is considered very stable, and should be useful for solving problems that are not well-posed, or that are well-posed but ill-conditioned, especially if stronger convergence criterion is needed.

Using this methodology, together with the rise in computing power, it is possible to establish exact or (with arbitrary precision arithmetic) very precise and more realistic benchmarks for algorithm performance testing.

**Acknowledgement**

This work was fully supported by the Croatian Science Foundation under the project IP–2014–09–2899.